\documentclass{elsart}

\usepackage{amssymb, color, epsfig}
\usepackage[all]{xy}
\input{psfig.sty}
\xyoption{dvips}

\newcommand{\abs}[1]{{\left\vert #1 \right\vert}}

\def\a{\alpha}

\def\g{\gamma}

\def\h{\eta}

\def\k{\kappa}

\def\x{\xi}

\def\io{\iota}

\def\ph{\phi}

\def\D{\Delta}

\def\Om{\Omega}

\def\addots{\mathinner{\mkern1mu\raise1pt\vbox{\kern7pt\hbox{.}}\mkern2mu    \raise4pt\hbox{.}\mkern2mu\raise7pt\hbox{.}\mkern1mu}} 

\def\({\left(}
\def\){\right)}

\def\del{\partial}

\def\bo#1{{\bf #1}}

\def\RR{\bo R}
\def\ZZ{\bo Z}

\newcommand{\id}{{\mathrm{id}}}
\newcommand{\emb}{{\mathrm{Emb}}}
\newcommand{\imm}{{\mathrm{Imm}}}
\newcommand{\hofiber}{{\mathrm{hofiber}}}
\newcommand{\tfiber}{{\mathrm{tfiber}}}
\newcommand{\holim}{{\mathrm{holim}}}

\newcommand{\ivmap}{{\mathrm{ivmap}}}

\newcommand{\equ}{{\mathrm{Equ}}}
\newcommand{\map}{{\mathrm{map}}}

\newcommand{\F}{{\bf{F}}}

\begin{document}

\begin{frontmatter}



\title{Embeddings in the 3/4 range}


\author{Brian A. Munson}
\ead{\texttt{munson@math.stanford.edu}}
\address{Department of Mathematics, Stanford University, Stanford, CA 94305, USA}

\begin{abstract}
We give a complete obstruction to turning an immersion $f: M^m\rightarrow\RR^{n}$ into an embedding when $3n\geq4m+5$. It is a secondary obstruction, and exists only when the primary obstruction, due to Andr\'{e} Haefliger, vanishes. The obstruction lives in a twisted cobordism group, and its vanishing implies the existence of an embedding in the regular homotopy class of $f$ in the range indicated. We use Tom Goodwillie's calculus of functors, following Michael Weiss, to help organize and prove the result.
\end{abstract}

\end{frontmatter}


\section{Introduction}

The story of immersions and embeddings of smooth manifolds begins with Whitney in 1936, when he proved his so-called ``easy'' embedding theorem: 
\begin{thm}\label{Whit}\cite{wh1}
Suppose $n\geq2m+1$. If $f:M^m\rightarrow N^n$, then $f$ is homotopic to an embedding.
\end{thm}


In \cite{wh2} he proved that every smooth manifold $M^m$ immerses in $\RR^{2m-1}$ and embeds in $\RR^{2m}$. There are obstructions in both cases: for immersions the proposed map might have singularities, and for embeddings the map might have self-intersections. In both cases he came up with a geometric elimination of the obstruction when it vanishes algebraically, and algebraic vanishing is automatic in this case. For maps into $\RR^{2m}$, this method of elimination is known as the Whitney trick. If we consider embeddings of $M^m$ in $N^{2m}$, there is a further obstruction to carrying out the Whitney trick if $N$ is not simply connected. More generally, we can use the same ideas to eliminate intersections of a $p$-manifold and $q$-manifold in a $(p+q)$-manifold. One application of this version of the Whitney trick is in the proof of the $h$-cobordism theorem, a corollary of which is the Poincar\'e conjecture in dimensions five and higher.  The reason for the dimensional restriction is that the Whitney trick works only for $p,q>2$.

In 1962, Haefliger \cite{haef} generalized the Whitney trick in a range of dimensions. Haefliger's assumption on dimension assures the immersion in question has no triple points, but the precise number really depends on making sure the Whitney trick will work. 

\begin{defn}
A map $F:M\times M\rightarrow N\times N$ is called \textsl{isovariant} if it is $\Sigma_2$-equivariant with respect to the action which switches the coordinates, and has the property that $F^{-1}(\Delta_N)=\Delta_M$ and $DF^{-1}(\Delta_{TN})=\Delta_{TM}$, where $\Delta_X$ denotes the diagonal of $X\times X$. We denote the space of isovariant maps by $\ivmap^{\Sigma_2}(M\times M,N\times N)$.
\end{defn}

For example, if $f$ is an embedding, then $f\times f$ is isovariant. To ease the statements of subsequent theorems and for the purposes of this paper we specialize to $N=\RR^n$.

\begin{thm}\label{Haef}\cite{haef}
Suppose $2n\geq3m+3$. Let $g:M^m\rightarrow \RR^{n}$ be an immersion, and suppose there exists an isovariant map $F(x,y): M^m\times M^m \rightarrow \RR^n\times \RR^n$ and an equivariant homotopy from $F$ to $g\times g$. Then $g$ is regularly homotopic to an embedding.
\end{thm}

There is a map $\emb(M,\RR^n)\rightarrow\ivmap^{\Sigma_2}(M\times M,\RR^n\times\RR^n)$ given by $f\mapsto f\times f$, and Theorem \ref{Haef} says that this map is $0$-connected. Haefliger also shows in \cite{haef} that if $2n>3m+3$, then the map $f\mapsto f\times f$ is $1$-connected. There is a further improvement due to J.-P. Dax \cite{dax}. It states

\begin{thm}\label{Dax}
The map $\emb(M,\RR^n)\rightarrow\ivmap^{\Sigma_2}(M\times M,\RR^n\times\RR^n)$ given by $f\mapsto f\times f$ is $(2n-3m-3)$-connected.
\end{thm}

Similar statements are true with a generic smooth manifold $N$ in place of $\RR^n$.







Dax's improvement is interesting because it gives a stable range description of the space of embeddings in terms of something more homotopy theoretic. This reduction of questions in differential topology to questions in homotopy theory is very much in the spirit of the Smale-Hirsch theorem for immersions. The approximation in question here replaces the global condition that an embedding $f$ should send a distinct pair of points to a distinct pair of points by a local property, that $F$ should take off-diagonal points to off-diagonal points.



Following \cite{we1}, we adopt a more modern view of spaces of embeddings where we will apply the calculus of functors. Denote by $\emb(M,N)$ the space of embeddings of $M^m$ in $N^n$, and the corresponding space of immersions by $\imm(M,N)$. We assume $m<n$. There is an inclusion $\emb(M,N)\rightarrow \imm(M,N)$. We are interested in the following types of questions: Given $g\in \imm(M,N)$, does it lift to $\emb(M,N)$? Is $\emb(M,N)$ non empty? What are the homotopy groups of $\emb(M,N)$, and how can we calculate them? The idea is to consider $\emb(-,N)$ as a \emph{cofunctor} (contravariant functor) from the poset $\mathcal{O}(M)$ of open subsets of $M$ to the category of spaces, $V\mapsto \emb(V,N)$. Theorems of Goodwillie, John Klein, and Michael Weiss (\cite{gkw2} and \cite{gw}) say that when $n-m>2$, there is a map from $\emb(M,N)$ to a space made from $\emb(V,N)$, where $V$ ranges over open subsets diffeomorphic to at most $k$ open balls, whose connectivity increases with $k$ (see below). We can understand embeddings of $k$ distinct balls in terms of configurations spaces of $k$ points plus some tangential information.

The \textsl{Taylor tower} of the embedding cofunctor is a sequence of cofunctors $\mathcal{T}_k\emb(V,N)$ with maps $\mathcal{T}_k\emb(V,N)\rightarrow \mathcal{T}_{k-1}\emb(V,N)$, where $V\in\mathcal{O}(M)$. We abbreviate $\mathcal{T}_k\emb(M,N)$ by $\mathcal{T}_k$. The spaces $\mathcal{T}_k$ are piecemeal descriptions of $\emb(M,N)$ in the sense that they only consider ``compatible'' embeddings of $k$ disjoint balls in $M$. A useful case to think about is $\mathcal{T}_1\emb(M,N)$, for it turns out to be homotopy equivalent to $\imm(M,N)$. We define

$$\mathcal{T}_1\emb(M,N)=\holim_{V\cong B^n}\emb(V,N).$$

Ignoring the homotopy inverse limit, observe that the inclusion $\emb(V,N)\rightarrow\imm(V,N)$ is a homotopy equivalence when $V$ is diffeomorphic to an open ball. So replacing $\emb(V,N)$ with $\imm(V,N)$ above, it remains to see that $\mathcal{T}_1\imm(M,N)$ is homotopy equivalent to $\imm(M,N)$. This fact is a reformulation of the Smale-Hirsch theorem. We say then that the first degree Taylor approximation to the space of embeddings is the space of immersions.

This Taylor approximation improves as $k$ gets large, provided that the codimension $n-m>2$. In fact, Goodwillie and Klein \cite{gk} prove that the map $\emb(M,N)\rightarrow \mathcal{T}_k\emb(M,N)$ is $(k(n-m-2)+1-m)$-connected. Taking $k=1$ and our note above about $\mathcal{T}_1$, we see that this says the map $\emb(M,N)\rightarrow\imm(M,N)$ is $(n-2m-1)$-connected, an improved version of Whitney's Theorem, which we stated as Theorem \ref{Whit}. From this setup we can also deduce Dax's  improvement of Haefliger's Theorem \ref{Haef}. If we take $k=2$, then the map $\emb(M,N)\rightarrow \mathcal{T}_2$ is $(2n-3m-3)$-connected. Goodwillie-Klein-Weiss \cite{gkw1} show that that $\mathcal{T}_2\emb(M,N)$ is equivalent to Haefliger's approximation (the space of isovariant maps, see Theorem \ref{Haef}) to the space of embeddings, which is of most interest when $2n\geq3m+3$. So when $2n\geq3m+3$, the problem of turning an immersion into an embedding is equivalent to studying the existence of liftings of elements of $\mathcal{T}_1$ to $\mathcal{T}_2$, liftings of immersions to isovariant maps. The next natural thing to consider is the case $k=3$, and the map $\emb\rightarrow\mathcal{T}_3$, which can produce embeddings when $3n\geq4m+5$ according to these connectivity estimates. Our Theorem \ref{Muns} concerns liftings from $\mathcal{T}_2$ to $\mathcal{T}_3$.

Before we state our Theorem \ref{Muns}, it will be useful to reformulate Dax's improvement (Theorem \ref{Dax}) of Haefliger's Theorem \ref{Haef} in terms of cubical diagrams and ``cobordism spaces'' (see Section \ref{prelim} for more information about cubical diagrams and Section \ref{cobmodel} for more details about cobordism spaces). Dax himself uses cobordism groups in \cite{dax}.

Haefliger's theorem \ref{Haef} tells us when the elimination of the double point obstruction is enough to produce an embedding. Given an immersion $g:M^m\rightarrow N^{n}$, consider $g\times g:M\times M\rightarrow N\times N$. We may assume that $g\times g$ is transverse to $\Delta_N$, and thus $(g\times g)^{-1}(\Delta_N)\setminus\Delta_M$ is a compact $(2m-n)$-dimensional submanifold of $M\times M\setminus\Delta_M$. The equivariant homotopy between $g\times g$ and the isovariant map $F$ can be regarded as a null-bordism of the double point set, because $F^{-1}(\Delta_N)\setminus\Delta_M=\emptyset$. Haefliger's Theorem \ref{Haef} says that when $2n\geq3m+3$, a null-bordism of the double point set is enough to produce an embedding in the regular homotopy class of $g$. We are now ready recast Dax's Theorem \ref{Dax} in terms of cobordism.

There is a simplicial set $C_2(M^m)$ which is a ``cobordism space'' in the sense that the homotopy groups of its realization are cobordism groups in which the double point obstruction lies. In this case, $\pi_k\abs{C_2(M)}\cong\Om_{2m-n+k}^{nL-T{M\choose2}}{M\choose 2}$ (see section \ref{cobmodel} for information about this notation. $M\choose k$ denotes the quotient by the $\Sigma_k$ action of $M^k\setminus\Delta$, where $\Delta$ is the fat diagonal). The map from the space of immersions $\mathcal{T}_1\emb(M^m,\RR^n)$ to $C_2(M)$ is defined by sending an immersion to its double point set, and the map $\ast\rightarrow C_(M)$ maps to the empty manifold.


\begin{thm}\label{haefrecast}
The following square is $(2n-3m-3)$-cartesian:

$$\xymatrix{
\emb(M^m,\RR^n) \ar[rr]\ar[dd] & &  \ast \ar[dd] \\
 & & \\
\mathcal{T}_1\emb(M^m,\RR^n)  \ar[rr]& & C_2(M^m)\\
}
$$
\end{thm}

In particular, if $2n\geq3m+3$, then the map $$\emb(M^m,\RR^n)\rightarrow\holim\left(\mathcal{T}_1\emb(M^m,\RR^n)\rightarrow C_2(M^m)\leftarrow\ast\right)$$ is onto $\pi_0$, and hence an immersion together with a cobordism to the empty manifold are enough to produce an embedding. More precisely,

\begin{cor}
Given an immersion map $f:M^m\rightarrow\RR^n$, there exists a manifold

$$D\in\Om_{2m-n}^{nL-T{{M}\choose2}}{{M}\choose2}$$

which represents the obstruction to lifting $g$ to $\mathcal{T}_2\emb(M^m,\RR^n)$. If $2n-3m-3\geq 0$, and $D$ is null-cobordant then there exists an embedding of $M$ in $\RR^n$.
\end{cor}

To construct embeddings of $M^m$ in $\RR^n$ in the range $3n\geq4m+5$, it is enough to produce an element of $\mathcal{T}_3$ to produce an embedding. We focus on lifting from $\mathcal{T}_2$ to $\mathcal{T}_3$, and use Haefliger's Theorem \ref{Haef} to interpret an element of $\mathcal{T}_2$ as an isovariant map $F:M\times M\rightarrow\RR^n$. The map $\k$ from $\mathcal{T}_2\emb(M^m,\RR^n)$ to a cobordism space $C_3(M)$ is the map which associates to each isovariant map $F$ and triple of points $(x_1,x_2,x_3)\in M^3\setminus\Delta$ the submanifold of $M^3\setminus\Delta$ where the three vectors $F(x_2,x_3), F(x_3,x_1)$ and $F(x_1,x_2)$ point in the same direction (compare \cite{bor}). This definition is less intuitive than that of the obstruction defined in \cite{mun}, where the construction of the obstruction class is obtained by following Haefliger's proof of Theorem \ref{Haef}, and we hope to have this written up soon for publication. The obstruction given in this paper has the advantage of being easier to define, and the computations of section \ref{generator} show that these two are equivalent in the sense that the two classes are cobordant. In fact, the way we discovered the definition of the obstruction presented here was to follow our work in \cite{mun} and guess manifolds of the right dimension until we found one that worked.

\begin{thm}\label{Muns}
The following square is $(3n-4m-5)$-cartesian:

$$\xymatrix{
\emb(M^m,\RR^n) \ar[rr]\ar[dd] & &  \ast \ar[dd] \\
 & & \\
\mathcal{T}_2\emb(M^m,\RR^n)  \ar[rr]^-{\k}& & C_3(M^m)\\
}
$$
\end{thm}

As above, an immediately corollary is

\begin{cor}
Given an isovariant map $F\in \mathcal{T}_2\emb(M^m,\RR^n)$, there exists a manifold

$$Z\in\Om_{3m-2n+2}^{(n-1)P-T{{M}\choose3}}{{M}\choose3}$$

which represents the obstruction to lifting $F$ to $\mathcal{T}_3\emb(M^m,\RR^n)$. If $3n-4m-5\geq 0$, then if $Z$ is null-cobordant there exists an embedding of $M$ in $\RR^n$.
\end{cor}

If the map $F$ is a lift of an immersion $g$, then if $Z$ is null-cobordant, there is an embedding in the regular homotopy class of $g$. There is no such embedding if and only if every lift of $g$ to $\mathcal{T}_2\emb(M^m,\RR^n)$ gives a non-trivial element of this group. An induction argument inspired by the proof of Theorem 5.1 in \cite{we1} reduces Theorem \ref{Muns} to the case where $M$ consists of exactly three points. The bulk of the proof is spent proving this special case, where we have to make some explicit calculations with the map $\k$. It is an instructive exercise to carry out a proof of Theorem \ref{haefrecast} in the same manner we prove our Theorem \ref{Muns}.


The obstruction $\tau$ defined in \cite{st} for immersions of a $2$-sphere in a $4$-manifold is the same as our obstruction $Z$ when the $4$-manifold in question is $\RR^4$. A generalization of our Theorem \ref{Muns} to embeddings in manifolds should make the connection between these two complete, which we are currently working on.


\subsection{Conventions}\label{conventions}

We write $QX$ for $\Om^{\infty}\Sigma^{\infty}X$ where X is a based space. We write $M^k\setminus\Delta$ for the complement of the fat diagonal in $M^k$. When we say a map is an \textsl{equivalence}, we mean it is a weak equivalence, unless otherwise noted. For a vector bundle $\xi$ over a space X, we denote by $T(\xi)$ is Thom space. We write Spaces for the category of fibrant simplicial sets, and we work in this category unless otherwise noted. A $k$-simplex in $\emb(M,N)$ is a fiber-preserving embedding of $M\times\Delta^k\rightarrow N\times\Delta^k$. By fiber-preserving we mean that if $f_k$ is a $k$-simplex of $\emb(M,N)$ and $p_N:N\times\Delta^k\rightarrow\Delta^k$ is the projection, then the composition $p_N\circ f=p_M$, where $p_M:M\times\Delta^k\rightarrow\Delta^k$ is the projection.



\section{Preliminary material}\label{prelim}

Our discussion of cubical diagrams is based on material from \cite{tg2}, and our discussion of the calculus of functors and spaces of embeddings is based on material from sections 0, 1, and 2 from \cite{we1}. The reader should look to these references for more details.

\subsection{Cubical diagrams}\label{cubicaldiagrams}

Cubical diagrams play a central role in the calculus of functors. We give the basic definitions and a brief discussion of their meaning.

\begin{defn}\label{ncube}
An \textsl{n-cube} of spaces is a functor $X$ from the category $\mathcal{P}_n$ of subsets of $\{1,\ldots , n\}$ to the category of spaces. We denote the value of $X$ at an object $S$ of $\mathcal{P}_n$ by $X_S$.
\end{defn}

Thus a $0$-cube is a space, a $1$-cube is a map of spaces, and a $2$-cube is a commutative square diagram.

\begin{defn}[1.3 of \cite{tg2}]\label{cartesian}
The $n$-cube $X$ is \textsl{homotopy cartesian} if the map $a(X):X_{\emptyset}\rightarrow\holim_{S\neq\emptyset}X_S$ is a weak equivalence. We say the cube is \textsl{$k$-cartesian} if the map $a(X)$ is $k$-connected map.
\end{defn}

\begin{defn}[1.1b of \cite{tg2}]
If $X$ is an $n$-cube of based spaces, we define the \textsl{total homotopy fiber} of $X$ as $\hofiber(a(X))$, and denote this space by $\tfiber(X)$.
\end{defn}

An immediate consequence of these last two definitions is that a cubical diagram $X$ is $k$-cartesian if and only if $\tfiber(X)$ is $(k-1)$-connected. One can also think of the total homotopy fiber as an inductive homotopy fiber. That is, view an $n$-cube $X$ as a map of $(n-1)$-cubes $Y\rightarrow Z$, and define $\tfiber(X)$ as $\hofiber(\tfiber(Y)\rightarrow\tfiber(Z))$. For a $0$-cube, define $\tfiber(X)=X$. See the beginning of Section 1 of \cite{tg2} for more details.


For example, Haefliger's Theorem \ref{Haef} states that the $2$-cube

$$\xymatrix{
\emb(M^m,\RR^n) \ar[rr]\ar[dd] & &  \ast \ar[dd] \\
 & & \\
T_1\emb(M^m,\RR^n)  \ar[rr]& & C_2(M^m)\\
}
$$

is $(2n-3m-3)$-cartesian. This means the map $$\emb(M^m,\RR^n)\rightarrow\holim(T_1\emb(M^m,\RR^n)\rightarrow C_2(M^m)\leftarrow \ast)$$

is $(2n-3m-3)$-connected. Recall that a point in $\holim(X\rightarrow Z\leftarrow Y)$ is a point in $X$, a point in $Y$, and a path between their images in $Z$. If $2n\geq 3m+3$, the above map is onto for $\pi_0$, and to produce an embedding it is enough to produce an immersion - an element of $T_1\emb(M^m,\RR^n)$ - whose double point manifold is null-cobordant. Equivalently we can interpret this theorem as saying that there is a $(2n-3m-3)$-connected map $$\hofiber\left(\emb(M^m,\RR^n)\rightarrow T_1\emb(M^m,\RR^n)\rightarrow \hofiber(\ast\rightarrow C_2(M^m))\right).$$ Since $\hofiber(\ast\rightarrow C_2(M^m)\simeq \Om C_2(M^m)$, this is another way of saying that the difference between embeddings and immersions in the range $2n\geq 3m+3$ is a double point obstruction.

\subsection{Calculus of functors and spaces of embeddings}

Let $M$ be a smooth manifold, and let $F:\mathcal{O}(M)\rightarrow \mbox{ Spaces }$ be a contravariant functor (which we refer to as a \emph{cofunctor}). 

\begin{defn}
Let $V_1$ and $V_2$ be smooth manifolds with boundary. A codimension $0$ embedding $i_1:V_1\rightarrow V_2$ is called an \textsl{isotopy equivalence} if there is a codimension $0$ embedding $i_2:V_2\rightarrow V_1$ such that $i_2\circ i_1$ and $i_1\circ i_2$ are isotopic to $\id_{V_1}$ and $\id_{V_2}$ respectively.
\end{defn}

\begin{defn}
A cofunctor $F:\mathcal{O}(M)\rightarrow \mbox{ Spaces }$ is called \textsl{good} if (a) it takes isotopy equivalences to homotopy equivalences, and (b) if $V_i\subset V_{i+1}$ is a sequence of objects then $F(\cup_i V_i)\rightarrow \holim_iF(V_i)$ is a homotopy equivalence.
\end{defn}


\begin{rem}\label{goodremark}
Proposition 1.4 of \cite{we1} says both $\emb(-,N)$ and $\imm(-,N)$ are good cofunctors. Part (b) in the definition of good guarantees that the values of $F$ are completely determined by its values on compact codimension $0$ handlebodies, because we may write any open $V$ as a union of $V_i$ such that $V_i\subset V_{i+1}$, $V_i$ is the interior of a compact codimesion $0$ handlebody, and $\cup_i V_i=V$. For the purposes of this paper, however, we are not interested in values of functors on generic open sets, but only on those open sets which are the interiors of smooth compact handlebodies. Therefore we will define the value of a cofunctor $F$ satisfying (a) on a generic open set $V$ by $F(V)=\holim_iF(V_i)$, where $V_i\subset V_{i+1}$, $V_i$ is the interior of a compact codimesion $0$ handlebody, and $\cup_i V_i=V$. Hence we will only check part (a) in the future.
\end{rem}


\begin{defn}
For a good cofunctor $F$ we define the \textsl{$k^{th}$ Taylor approximation} to $F$, denoted $\mathcal{T}_kF:\mathcal{O}(M)\rightarrow\mbox{ Spaces }$, by

$$\mathcal{T}_kF(U)=\holim_{V\in \mathcal{O}_k(U)}F(V).$$

Here $\mathcal{O}_k(U)$ is the subcategory of $\mathcal{O}(U)$ consisting of those open sets $V\subset U$ which are diffeomorphic to at most $k$ open balls.
\end{defn}

\begin{defn}
We say that $F$ is \emph{polynomial of degree $\leq k$} if given pairwise disjoint closed subsets $A_0,A_1,\ldots,A_k$ of $U\in\mathcal{O}(M)$, the $(k+1)$-cube

$$S\mapsto F(U\setminus \cup_{i\in S}A_i)$$

is homotopy cartesian, where $S$ ranges through subsets of $\{0,1,\ldots,k\}$.
\end{defn}

The next two theorems state that the functors $\mathcal{T}_kF$ are polynomial and that they are essentially determined by their values on special open sets. 

\begin{thm}[\cite{we1}, Theorem 6.1]\label{tkpoly}
The cofunctor $\mathcal{T}_kF$ is polynomial of degree $\leq k$.
\end{thm}

\begin{thm}[\cite{we1}, Theorem 5.1]\label{polydet}
Suppose that $\g:F_1\rightarrow F_2$ is a morphism of good cofunctors, and that $F_i$ is polynomial of degree $k$ for $i=1,2$. Then if $\g:F_1(V)\rightarrow F_2(V)$ is a homotopy equivalence for all $V\in \mathcal{O}_k(M)$, then it is a homotopy equivalence for all $V\in \mathcal{O}(M)$.
\end{thm}

From its definition we see that the values of $\mathcal{T}_kF$ are completely determined by its values on $\mathcal{O}_k(M)$, so Theorem \ref{polydet} is not too surprising. The proof of this theorem inspired that of Theorem \ref{Muns}.

\subsection{A model for $\mathcal{T}_2\emb(M,\RR^n)$}\label{ttwomodel}


In \cite{gkw1}, the authors show that the homotopy pullback of 

$$\xymatrix{
  & \ivmap^{\Sigma_2}(M\times M, N\times N) \ar[d]\\
\map(M,N) \ar[r]_-{f\mapsto f\times f} &\map^{\Sigma_2}(M\times M, N\times N)\\
}
$$

is homotopy equivalent to $\mathcal{T}_2\emb(M,N)$. In the case $N=\RR^n$, the bottom two spaces are contractible, and thus

$$\mathcal{T}_2\emb(M,\RR^n)\simeq \ivmap^{\Sigma_2}(M\times M, \RR^n\times \RR^n).$$

We go further and replace $\ivmap^{\Sigma_2}(M\times M, \RR^n\times \RR^n)$ by the homotopy equivalent space $\ivmap^{\Sigma_2}(M\times M, \RR^n)$, where the $\Sigma_2$ action on $\RR^n$ is given by the antipodal map. The homotopy equivalence is given by the map $(f_1,f_2)\mapsto f_1-f_2$, with homotopy inverse $f\mapsto(f/2, -f/2)$ (a straight line homotopy will suffice here). The simplicial structure here is similar to that for the embedding space: a $k$-simplex $F^k$ of $\ivmap^{\Sigma_2}(M\times M, \RR^n)$ is a fiber-preserving isovariant map $M\times M\times\Delta^k\rightarrow \RR^n\times\Delta^k$, where $\Sigma_2$ acts trivially on $\Delta^k$, and by isovariance we mean that ${F^{k}}^{-1}(0\times\Delta^k)=\Delta_M\times\Delta^k$. The map from $\mathcal{T}_2$ to $\mathcal{T}_1$ is the map which restricts an isovariant map $F:M\times M\rightarrow \RR^n$ to the diagonal and records the induced map of normal bundles, which we may interpret as $TM$ and $T\RR^n$ respectively. Then using the Smale-Hirsch theorem, this gives an element of $\mathcal{T}_1$.





\section{Cobordism spaces}\label{cobmodel}

Since we have opted for a cobordism description of our obstruction, it will be useful to consider a cofunctor $C:\mathcal{O}(M)\rightarrow\mbox{ Spaces}$, which gives us a ``cobordism space'': a simplicial set whose realization has as its homotopy groups the cobordism groups we encounter in defining our obstruction. These groups will be denoted $\Om_{d+k}^{\x-\h}(X)$, where $X$ is a space and $\x$ and $\h$ are bundles on $X$. A zero simplex in this space is a triple $(W^d,f,\ph)$ (sometimes denoted by just $W$) where $W$ is a $d$-dimensional smooth manifold, $f:W \rightarrow X$ is continuous and proper, and $\ph$ is a stable isomorphism $TW \oplus f^*\x \cong f^*\h$. The equivalence relation is the usual one defined by $(d+1)$-dimensional manifolds with boundary. The equivalence class of the $W$ described above defines an element of the group $\Om_{d}^{\x-\h}(X)$. There is an equivariant version of these spaces and groups, which we pause to mention because we use it in the proof of Theorem \ref{Muns}.


Let $\widetilde{X}$ be a space with a free $G$ action for some group $G$, and let $X$ denote the quotient by that $G$ action. Let $\widetilde{\xi}$ and $\widetilde{\h}$ be vector bundles on $\widetilde{X}$ with a $G$ action, and let $\xi$ and $\h$ be the quotient bundles on $X$. By abuse of notation, omitting $G$, a zero simplex of $C_{d}^{\xi-\h}(X)$ is a smooth, closed, compact manifold $W^d$ with continuous proper map $f:W\rightarrow X$ and a stable isomorphism $TW\oplus f^*(\xi)\cong f^*(\h)$. An element of $\Om_{d}^{\xi-\h}(X)$ is an equivalence class of such data.

Now suppose we are given a smooth closed manifold $W^d$ with free $G$ action, a continuous $G$-map $f:W\rightarrow \widetilde{X}$, and a stable $G$-isomorphism $\Phi:TW\oplus f^*(\widetilde{\xi})\rightarrow f^*(\widetilde{\h})$. The manifold $W/G$ is a zero simplex in $C_{d}^{\xi-\h}(X)$. More generally suppose that $H$ is a subgroup of $G$, and that $W$ is as above, only now $W$ has free $H$ action, $f$ is an $H$-map, and $\Phi$ is a stable $H$-isomorphism. Then $(G\times_H W)/G$ represents a zero simplex of $C_{d}^{\xi-\h}(X)$. We identify $G\times_H W$ with $\abs{G}/\abs{H}$ disjoint copies of the same manifold, now made into a $G$-space, with $G$-maps induced by the given $H$-maps. We are going to construct a cobordism class with $G=\Sigma_3$ and $H$ as one of the three copies of $\Sigma_2$.
 
From the description above we make a simplicial model $C_{\bullet}^{\xi-\h}(X)$ for a space $C_{d}^{\xi-\h}(X)$ whose realization has as its homotopy groups the cobordism groups mentioned above; $\pi_k\abs{C_{d}^{\xi-\h}(X)}=\Om_{d+k}^{\xi-\h}(X)$. It is related to the Thom space of a virtual bundle; see the remark following Proposition \ref{CKan}. Although this notation expresses the dependence on $d,\xi,\h$ and $X$, it is rather cumbersome, so we will usually omit it and just name the relevant parameters once and for all.

\begin{defn}[Simplicial Model for a Cobordism Space]\label{cobspacedefn} The simplicial set $C_{\bullet}^{\xi-\h}(X)$ has as its $0$-simplices the set $C_0=\{(W^d,f,\ph)\}$, where $W$ is embedded in $\RR^{\infty}$, $f:W\rightarrow X$ is a continuous and proper map, and $\ph$ is a stable isomorphism $\ph: TW\oplus f^*(\xi)\rightarrow f^*(\eta)$. The $1$-simplices are
$C_1=\{(W^{d+1},f,\ph)\}$ where $W$ is embedded in $\RR^{\infty}\times \D^1$, $W$ is transverse to $\RR^{\infty}\times \del\D^1$, $f:W\rightarrow X$ is continuous and proper, and $\ph: TW\oplus f^*(\xi)\rightarrow f^*(\eta)$ is a stable isomorphism. In general, the $k$-simplices are the set $C_k=\{(W^{d+k},f,\ph)\}$ where $W$ is embedded in $\RR^{\infty}\times \D^k$, $W$ is transverse to $\RR^{\infty}\times \del_S\D^k$ for all nonempty subsets $S\subset\{0,1,\ldots, k\}$, $f:W\rightarrow X$ is continuous and proper, and $\ph: TW\oplus f^*(\xi)\rightarrow f^*(\eta)$ is a stable isomorphism.
\end{defn}

We will also make use of a relative version of this construction for a pair $(X,Y)$. A $k$-simplex of $C_{\bullet}^{\xi-\h}(X,Y)$ is a $(k+d)$-dimensional manifold $W$ with boundary $\del W$ embedded in $\RR^{\infty}\times\Delta^k$ such that $W$ and $\del W$ are transverse to $\del_S\Delta^k$ for all nonempty $S\subset\{1,\ldots, k\}$. Moreover, there is a continuous proper map of pairs $f:(W,\del W)\rightarrow (X,Y)$, and stable isomorphisms $TW\oplus f^*(\xi)\rightarrow f^*(\eta)$ and $T\del W\oplus f^*(\xi)\rightarrow f^*(\eta)$ which are compatible in the sense that there is a commutative diagram relating the bundle isomorphisms on $W$ and $\del W$. Now a $k$-simplex $W$ has boundary $\del W$, and the boundary defines a $(k-1)$-simplex of $C_\bullet(Y)$. In the case $Y=\emptyset$, $C_\bullet(X,Y)=C_\bullet(X)$ (which forces $\del W=\emptyset$).


Moreover, the manifolds $W^{d+k}\subset \Delta^k\times\RR^\infty$ should the \emph{conditioned}. To be conditioned means that if we denote by $W_t$ the part of $W$ that sits over $t\in\Delta^k$, then $W_t$ should be independent of $t$ in a neighborhood of $\cup_i \del_i\Delta^k$.

The face and degeneracy maps are induced by those of $\D^{\bullet}$. The $i^{\mbox{th}}$ face map $d_i:C_k\rightarrow C_{k-1}$ is just the intersection of $W^{d+k}$ with the $i^{\mbox{th}}$ face of $\D^k$. The $i^{\mbox{th}}$ degeneracy map $s_i:C_k \rightarrow C_{k+1}$ takes $W$ to the fiber product $W'$

$$\xymatrix{
W' \ar[r]\ar[d] & \RR^{\infty}\times \D^{k+1}\ar[d]^{s_i} \\
W\ar[r] & \RR^{\infty}\times \D^{k}\\
}
$$

where $s_i$ is the $i^{th}$ degeneracy for $\D^{\bullet}$. That it satisfies the axioms for a simplicial set is straightforward, because we are building on the usual simplicial structure on $\Delta^k$.

\begin{prop}\label{CKan}
$C_{\bullet}^{\xi-\h}(X)$ is a Kan complex.
\end{prop}


\begin{pf}
Recall that a simplicial set $C_\bullet$ satisfies the Kan extension condition if for every collection of $k+1$ $k$-simplices $x_0, x_1,\ldots, x_{n-1},x_{n+1},\ldots, x_{k+1}$ satisfying $\del_ix_j=\del_{j-1}x_i$ for $i<j$, $i,j\neq n$, there exists a $(k+1)$-simplex $x$ such that $\del_ix=x_i$ for all $i\neq n$. Let $\Delta^{k+1}$ be embedded in $\RR^{k+1}$ in the standard way, and denote by $\del\Delta^{k+1}_{\widehat{n}}$ the union of all but the $n$th face $\del_n\Delta^{k+1}$ of $\Delta^{k+1}$. Let $r:\Delta^{k+1}\rightarrow \del\Delta^{k+1}_{\widehat{n}}$ be defined by $r(x)=y$ if $x$ is on the line perpendicular to $\del_n\Delta^{k+1}$ passing through $y$. It is well-defined because the restriction $p_{\widehat{n}}$ to $\del\Delta^{k+1}_{\widehat{n}}$ of the orthogonal projection $p$ onto the $k$-plane containing $\del_n\Delta^{k+1}$ is one-to-one. Let $W_0, W_1, \ldots, W_{n-1}, W_{n+1}, \ldots, W_{k+1}$ be a collection of $(k+1)$ $k$ simplices satisfying the hypotheses of the Kan extension condition. Define 

$$\widehat{W}=\bigcup_{\del_iW_j=\del_{j-1}W_i}W_i\subset \RR^{\infty}\times\del\Delta^{k+1}_{\widehat{n}}$$ 

The manifold $\widehat{W}$ defines a $k$-simplex itself if one identifies $\del\Delta^{k+1}_{\widehat{n}}$ with $\del_n\Delta^{k+1}$ using $p_{\widehat{n}}$. The map $\widehat{f}:\widehat{W}\rightarrow X$ is made by gluing together the $f_i:W_i\rightarrow X$ according to $\del_iW_j=\del_{j-1}W_i$. The map $\widehat{f}$ is proper because the $f_i$ are. The stable bundle isomorphism $T\widehat{W}\oplus \widehat{f}^\ast(\xi)\rightarrow \widehat{f}^\ast(\eta)$ is made in exactly the same way. Define $W$ by the fiber product

$$\xymatrix{
W \ar[r]\ar[d] & \RR^{\infty}\times \Delta^{k+1}\ar[d]^{\id\times r} \\
\widehat{W}\ar[r] &\RR^{\infty}\times\del\Delta^{k+1}_{\widehat{n}}\\
}
$$

Then $W$ defines a $(k+1)$-simplex.
$\Box$\end{pf}

\begin{rem}\label{cobspacethomspace}
This space is equivalent to $QT(\x-\h)$. For details on how to make sense of the Thom space of a virtual bundle, see \cite{gkw2}. To see this, consider the subcomplex of the total singular complex of $QT(\x-\h)$ consisting of those $k$-simplices $\k:\Delta^k\rightarrow \Om^n\Sigma^n(T(\x-\h))$ that correspond to maps $\k':\Sigma^n(\Delta^k)\rightarrow \Sigma^n(T(\x-\h))$ which are transverse to the zero section of $T(\x-\h)$. This sub-complex is equivalent to the full complex and the map $\k\mapsto \k'^{-1}(0)$ to the cobordism model is an equivalence. See \cite{tg1} for a similar construction. 
\end{rem}

That $C_\bullet$ is a Kan complex ensures that the homotopy groups of its realization will be the cobordism groups we want. It is also used to prove two of the next three propositions.

\begin{prop}
There is an equivalence

$$C_{d+l}^{\xi-\h}(X)\simeq \Om^lC_{d}^{\xi-\h}(X).$$
\end{prop}







\begin{pf}
We prove this in the case $l=1$, iterating to obtain the general case. We need the relative version of $C$ mentioned after Definition \ref{cobspacedefn}. There is a map $C_\bullet(X,Y)\rightarrow C_{\bullet -1}(Y)$ given by taking the boundary, and this map is a Kan fibration. One can adapt the proof of Proposition \ref{CKan} to the relative setting to check this, as a simplicial set is a Kan complex if an only if the map of it to a one-point complex is a Kan fibration. Since $C_{\bullet -1}(Y)$ is a Kan complex by Proposition \ref{CKan}, so then is $C(X,Y)$ a Kan complex, and its homotopy groups are the relative bordism groups. The fiber of this map is $C_\bullet(X)$, because this is precisely what maps to the basepoint in $C_{\bullet -1}(Y)$. Furthermore, since this map is a Kan fibration, $C_\bullet(X)$ is also equivalent to the homotopy fiber. If we specialize to the case $X=Y$, we have $C_\bullet(X)=\hofiber(C_\bullet(X,X)\rightarrow C_{\bullet-1}(X))$. Finally, observe that $C_\bullet(X,X)$ is contractible.
$\Box$\end{pf}


Now let us consider the special case when $X$ is a smooth manifold of dimension $k$. In this case we can consider $C$ as a cofunctor $C:\mathcal{O}(X)\rightarrow\mbox{ Spaces}$.

\begin{prop}\label{Cisgood}
$C:\mathcal{O}(X)\rightarrow\mbox{ Spaces}$ is a good cofunctor.
\end{prop}

\begin{pf}
We need to check that given open sets $U_1,U_2\in \mathcal{O}(X)$, with $U_1\subset U_2$, we get a map $C(U_2)\rightarrow C(U_1)$. Suppose then that we have a smooth manifold $M^k$ with a continuous proper map $f:M\rightarrow U_2$ with bundle data. We may assume that $f$ is smooth and transverse to $U_1\subset U_2$. Then $f^{-1}(U_1)$ is a smooth manifold of dimension $k$ and $f:f^{-1}(U_1)\rightarrow U_1$ is proper and $f^{-1}(U_1)$ has the right kind of bundle data too, since the bundle data it receives is that of $M$ pulled back to $f^{-1}(U_1)$.

To check part (a) of goodness, one can use exactly the reasoning Weiss uses for Proposition 1.4 in \cite{we1} applied to the functor $C$, and we refer the reader to Remark \ref{goodremark} for part (b).
$\Box$\end{pf}


\begin{prop}\label{cobtn}
Let $U\in\mathcal{O}(X)$ be a tubular neighborhood of a compact submanifold $S\subset X$, so that $U$ is a $k$-disk bundle over $S$. Then there is an equivalence $C_\bullet(U)\rightarrow C_{\bullet-k}(S)$, where we replace the bundle $f^*(\xi)$ by $f^*(\xi\oplus\nu(S\subset U))$ in the definition of $C_{\bullet-k}(S)$.
\end{prop}

\begin{pf}
Consider the sub-simplicial set $C_\bullet'(U)\subset C_\bullet(U)$ for which the map $W^{d+k}\rightarrow U$ is transverse to $S$. This subcomplex $C_\bullet'(U)$ is equivalent to $C_\bullet(U)$ (see Hypothesis 3.18 of \cite{tg1}). There is a map $i:C'_\bullet(U)\rightarrow C_{\bullet-k}(S)$ given by intersection with $S$. A $(d+k)$-simplex $W\in C_\bullet(U)$ gives a $d$-simplex $W\cap S\in C_{\bullet-k}(S)$ because the intersection is transverse. Moreover, there is a map $r:U\rightarrow S$ given by identifying $U$ as a tubular neighborhood of $S$ and sending $(s,v)\in U$ to $s\in S$. This induces a map $C(S)\rightarrow C(U)$ in the other direction. We claim that they are homotopy inverses. First, the composition $r\circ i:S\rightarrow U\rightarrow S$ is the identity. Given a $d$-simplex of $M$ of $S$, consider the fiber product $M\times_S U$. The mapping $r:U\rightarrow S$ is smooth, and we may assume the map $M\rightarrow S$ is smooth and transverse to $r$, so that $M\times_S U$ is a manifold with proper map to $U$. Since the composition above is the identity, the fiber product $(M\times_S U)\times_U S$ is equivalent to $M$, again with transversality assumptions. This process leaves the bundle data on $M$ alone in the sense that if $p_X$ denotes the canonical map from a pullback $X\times_Y Z$ to $X$, then $p_M^{\ast}p_{M\times_S U}^{\ast}$ is an isomorphism.

Now consider the composition $i\circ r:U\rightarrow S\rightarrow U.$ We have $(i\circ r)(s,v)=(s,0)$, and there is a homotopy $h:U\times I\rightarrow U$ from $i\circ r$ to $\id$ given by fiberwise retraction to the origin. Given a $(k+d)$-simplex $M$ of $C(U)$, the fiber product $M\times_U S$ is a $d$-simplex of $C(S)$ whose map to $S$ is proper, and again pulling back we get a $(k+d)$-simplex $M'=(M\times_U S)\times_S U$, again with proper map to $U$. The homotopy will provide us with a cobordism between $M$ and $M'$, as follows. Consider the fiber product

$$\xymatrix{
W \ar[r]\ar[d] & M\ar[d] \\
U\times I\ar[r]_{h} &U\\
}
$$

Again transversality assumptions ensure $W$ is a manifold. Since the map $M\rightarrow U$ is proper, so is the map $W\rightarrow U\times I$. If we denote by $W_t$ the submanifold of $W$ that sits over $U\times\{t\}$, then $W$ is a cobordism between $W_0=M'$ and $W_1=M$. The bundle data is pulled back in each step, and it is straightforward to check that this is the right bundle data in each case.









$\Box$\end{pf}

Let $X$ be a smooth manifold, and consider the space $C_{d}^{\xi-\h}{X\choose 3}$. Recall that ${X\choose 3}$ is the quotient by $\Sigma_3$ of $X^3\setminus\Delta$. We can also view this as a cofunctor $C:\mathcal{O}(X)\rightarrow\mbox{ Spaces}$, using the map $\mathcal{O}(X)\rightarrow\mathcal{O}({X\choose 3})$, although it is a bit awkward with this notation. Using our shorthand, for $U\in\mathcal{O}(X)$, we write $C_3(U)=C_{d}^{\xi-\h}{U\choose 3}$.

\begin{prop}\label{cobpoly3}
The cofunctor $C:\mathcal{O}(X)\rightarrow\mbox{ Spaces}$ defined by $C_3(U)=C_{d}^{\xi-\h}{U\choose 3}$ is a polynomial of degree $\leq 3$.
\end{prop}

\begin{pf}
By abuse of notation use the letter $C_3$ also for the realization of the simplicial set, and for brevity we abbreviate $C_3=C$. We need to prove for all $U\in \mathcal{O}(X)$ and for all pairwise disjoint subsets $A_0,A_1,A_2,A_3$ of $U$ that the 4-cube $S\mapsto C(U\setminus \cup_{i\in S}A_i)$ is homotopy cartesian, where $S\subset\{0,1,2,3\}$. Observe that ${U\choose 3}=\bigcup_i{U\setminus A_i\choose 3},$ because the sets $A_i$ are four in number. For every subset $S\subset\{0,1,2,3\}$, let $V_S={U\setminus \cup_{i\in S}A_i\choose 3}.$

Our goal is to show that $C(V_\emptyset)\rightarrow\holim_{S\neq \emptyset}C(V_S)$ is a homotopy equivalence. If $U_1$ and $U_2$ are open sets, then 

\begin{equation}\label{cobexcision}
\xymatrix
{
C(U_1\cup U_2) \ar[r]\ar[d] & C(U_1)\ar[d] \\
C(U_2)\ar[r] &C(U_1\cap U_2)\\
}
\end{equation}

is homotopy cartesian. This is a restatement of the fact that this cobordism cofunctor satisfies the excision axiom.

We view the $4$-cube $S\mapsto C(V_S)$ as a map of $3$-cubes. If we let $T$ range through subsets of $\{0,1,2\}$, then the map of $3$-cubes we have in mind is

$$(T\mapsto C(V_T))\rightarrow(T\cup\{3\}\mapsto C(V_{T\cup\{3\}})).$$

By proposition 1.6 of \cite{tg2}, it is enough to show that each of these $3$-cubes is homotopy cartesian to show that the entire $4$-cube is. The argument for both is exactly the same, so let us only indicate why this is true for $T\mapsto C(V_T)$. Since $T$ is ranging though subsets of $\{0,1,2\}$, we represent it as


\begin{equation}\label{cubicalcob1}
\xymatrix@=20pt
{
C(V_\emptyset)\ar[rr]\ar[dd]\ar[rd]         &           &   C(V_0) \ar'[d][dd]           \ar[dr]  &                  \\
        &  C(V_1) \ar[rr] \ar[dd]  &             & C(V_{01})
         \ar[dd] \\
C(V_2) \ar'[r][rr] \ar[dr] &        &   C(V_{02})
\ar[dr] &                   \\
        &    C(V_{12}) \ar[rr]      &                    &
 C(V_{012})
}
\end{equation}


Also consider the related diagram

\begin{equation}\label{cubicalcob2}
\xymatrix@=20pt
{
C(V_{\emptyset}) \ar[r]\ar[dd] & C(V_0\cup V_1)\ar[rr]\ar[dd]\ar[rd]         &           &   C(V_0) \ar'[d][dd]           \ar[dr]  &                  \\
&        &  C(V_1) \ar[rr] \ar[dd]  &             & C(V_{01})
         \ar[dd] \\
C(V_2) \ar[r] & C((V_0\cup V_1)\cap V_2) \ar'[r][rr] \ar[dr] &        &   C(V_{02})
\ar[dr] &                   \\
     &    &    C(V_{12}) \ar[rr]      &                    &
 C(V_{012})
}
\end{equation}


We wish to show that (\ref{cubicalcob1}) is homotopy cartesian. Since the $A_i$ are pairwise disjoint, $V_{ij}=V_i\cap V_j$, and hence each of these square faces of the cubical part of (\ref{cubicalcob2}) are homotopy cartesian, as they are special cases of (\ref{cobexcision}). Using Proposition 1.6 of \cite{tg2}, this proves that the cubical part of (\ref{cubicalcob2}) is homotopy cartesian. Notice that the square part of (\ref{cubicalcob2}) is homotopy cartesian because it is of the same form as (\ref{cobexcision}). Since both the cubical and square parts of (\ref{cubicalcob2}) are homotopy cartesian, it follows again from Proposition 1.6 of \cite{tg2} that (\ref{cubicalcob1}) is homotopy cartesian.
$\Box$\end{pf}

\subsection{Counting $0$-dimensional cobordism classes}

In section \ref{smoothknot} we need to identify the group $\Om_0^{\xi-\eta}(X)$. Suppose that $X$ is path connected and let $\pi=\pi_1(X)$. Let $x\in X$ be the basepoint of $X$, and let $\xi_x$ and $\eta_x$ denote the fibers of the bundles $\xi$ and $\eta$ over $x$. Give an orientation to $\xi_x$ and $\eta_x$ and let $\g\in\pi$. If we drag this orientation around the loop $\g$, we can ask whether or not the orientation class changed. This defines homomorphisms $w(\xi),w(\eta):\pi\rightarrow\{+1,-1\}$ which are $+1$ if the orientation class does not change, and $-1$ if it does. Note that choosing an orientation of $\xi_x$ is the same as choosing an isomorphism of $\xi_x$ with a trivial bundle over a point, up to homotopy.

\begin{prop}\label{countclasses}
The group $\Om_0^{\xi-\eta}(X)$ is isomorphic with $\ZZ$ if $w(\xi)=w(\eta)$, and $\ZZ/2$ if $w(\xi)\neq w(\eta)$.
\end{prop}

\begin{pf}
An element of $\Om_0^{\xi-\eta}(X)$ is represented by a finite set $S$ mapped to $X$ together with a stable isomorphism $\xi\cong_s\eta$ over $S$. A single point with necessarily trivial bundle data generates this group, and we may assume that this  point maps to the basepoint $x\in X$. Since both $\xi$ and $\eta$ become trivial over a point, there are two possible stable isomorphisms between them, classified by the sign of their determinants. Denote the two possible cobordism classes of a point by $+x=(x,f,\phi_+)$ and $-x=(x,f,\phi_-)$, where $\phi_+$ has positive determinant and $\phi_-$ has negative determinant, and $f$ is the inclusion of the basepoint $x$ in $X$. Both $+x$ and $-x$ represent generators of $\Om_0^{\xi-\eta}(X)$, and the proposition will follow when we show that $+x$ and $-x$ are cobordant if and only if $w(\xi)\neq w(\eta)$. Let $(I,F,\Phi)$ be a cobordism between $+x$ and $-x$. That is, $F:I\rightarrow X$ satisfies $F(0)=F(1)=x$, and $\Phi$ is an isomorphism $\Phi:TI\oplus F^*\xi\cong_s F^*\eta$, where $\cong_s$ denotes that the isomorphism is stable. We regard $\Phi$ as a homotopy, over $I$, between $\phi_+$ and $\phi_-$, and the only way they can have determinants of opposite sign is if exactly one of $\xi_x$ or $\eta_x$ had its orientation class change when dragged along $F$. If we let $g\in\pi$ denote the class defined by $F$, then $w(\xi)(g)\neq w(\eta)(g)$. Conversely, if $w(\xi)\neq w(\eta)$, then let $g\in\pi$ satisfy $w(\xi)(g)\neq w(\eta)(g)$. Choose a representative $\g:I\rightarrow X$ for $g$, where $\g(0)=\g(1)=x$. Then $\g$ gives rise to a cobordism between $+x$ and $-x$ as follows. The map $\phi_+$ can be interpreted as a choice of orientation of $\xi_x-\eta_x$. Dragging this orientation along $\g$ leads to a cobordism $(I,\g, \Phi)$ where the restriction of $\Phi$ to $0$ is $\phi_+$, and the restriction of $\Phi$ to $1$ is an isomorphism $\phi_-$ of negative determinant, since $w(\xi)(g)\neq w(\eta)(g)$. Hence $(I,\g,\Phi)$ is a cobordism between $+x$ and $-x$.
$\Box$\end{pf}

\subsection{The plane bundle $P$ and the cobordism space}\label{pbundle}

We now describe the specific cobordism space which arises in the statement of Theorem \ref{Muns}. Consider the trivial bundle $(M^3\setminus\Delta)\times\RR^2$ over $M^3\setminus\Delta$ with fibers $\RR^2$. Let $e_1,e_2,e_3$ be nonzero vectors in $\RR^2$ such that $e_1+e_2+e_3=0$. Let $\Sigma_3$ act linearly on $\RR^2$ by permuting these vectors. The quotient of this product by the $\Sigma_3$ action is the bundle $P$, which is a bundle over $M\choose 3$. Denote by $kP$ the $k$-fold direct sum of $P$. We let $\widetilde{P}$ denote the trivial $\RR^2$ bundle. We should mention that the line bundle $L$ over $M\choose 2$ mentioned in the introduction is made in an analogous way from the trivial rank 1 bundle on $M\times M\setminus \Delta$ by letting $\Sigma_2$ act by -1 on the fibers.

In Section \ref{cobmodel} we described a simplicial set $C_{d}^{\xi-\h}(X)$ such that $\pi_k\abs{C_{d}^{\xi-\h}(X)}=\Om_{d+k}^{\xi-\h}(X)$. In the case where $X$ is a smooth manifold, Proposition \ref{Cisgood} tells us that we may regard this space as one value of a good cofunctor $C:\mathcal{O}(X)\rightarrow\mbox{ Spaces}$.

\begin{defn}
Let $M^m$ be a smooth manifold. We define a cofunctor $C_3:\mathcal{O}(M)\rightarrow\mbox{ Spaces}$ by the rule $U\mapsto C_{3m-2n+2}^{(n-1)P-T{M\choose 3}}\left({U\choose 3}\right)$.
\end{defn}

\subsection{The map $\k:\mathcal{T}_2\emb(M^m,\RR^n)\rightarrow C_3(M^m)$}\label{themapkappa}

Recall that $\mathcal{T}_2\emb(M^m,\RR^n)\simeq\ivmap^{\Sigma_2}(M\times M,\RR^n)$. One feature of an element $F\in \ivmap^{\Sigma_2}(M\times M,\RR^n)$ is that for each pair of distinct points in $M$ it gives a nonzero vector in $\RR^n$. The map $\k:\mathcal{T}_2\emb(M^m,\RR^n)\rightarrow C_3(M)$ associates to each triple of distinct points in $M$ the submanifold of ${M\choose 3}$ where the three nonzero vectors determined $F$ point the same direction.

We begin by describing the map $\k$ for $0$-simplices. Consider the standard action of $\Sigma_3$ on the set $\{1,2,3\}$. Denote by $\Sigma_2^{ij}$ by the subgroup which switches $i$ and $j$ for $i\neq j$. 

\begin{defn}\label{aij}
Let $\RR_{>0}^3$ denote the open octant of $\RR^3$ where all three coordinates are positive. Denote points in this space by triples $\{(a_{23},a_{31},a_{12})\}$ with the $\Sigma_3$-action induced by its action on indices, where $a_{ij}=a_{ji}$.
\end{defn}

First consider the map $\F':(M^3\setminus\Delta)\times \RR_{>0}^2\rightarrow \RR^n\times\RR^n$ defined by 

$$
\F'(x_1,x_2,x_3,a_{12},a_{31})=(F(x_2,x_3)-a_{31}F(x_3,x_1),F(x_2,x_3)-a_{12}F(x_1,x_2))
$$

The zeros of this function occur when the $F(x_i,x_j)$ all point the same way since the $a_{ij}$ are all positive. To make the symmetric group action easier to analyze, we modify this map slightly.

\begin{defn}\label{Fdef}
Define $\F=(f+g,f-g)$, where we let $\F'=(f,g)$ be the map above.
\end{defn}

The map $\F$ is $\Sigma_2^{23}$-equivariant, where $\Sigma_2^{23}$ acts on $\RR^n\times \RR^n$ by $-1$ on the first factor and trivially on the second factor, and we may assume it is transverse to $0\times 0\in\RR^n\times\RR^n$ because the action of $\Sigma_2^{23}$ on $(M^3\setminus\Delta)\times \RR_{>0}^2$ is free.


\begin{defn}
Define $Z_1=\F^{-1}(0\times0)$.
\end{defn}

By transversality, $Z_1$ is a $3m-2n+2$-dimensional submanifold of $(M^3\setminus\Delta)\times \RR_{>0}^2$.

\begin{lem}\label{zbundle}
$Z_1$ is a compact, closed $3m-2n+2$-dimensional manifold with $\Sigma_2^{23}$ action and a $\Sigma_2^{23}$-equivariant map $p$ to $M^3\setminus\Delta$. Moreover, there is a $\Sigma_2^{23}$-equivariant isomorphism $TZ_1\oplus p^*nP\rightarrow p^*(T(M^3\setminus\Delta)\oplus P)$.
\end{lem}

\begin{pf}
The comments in the paragraph above give everything we need save compactness, that it is closed, and the bundle isomorphism. That it is closed follows from compactness since it is a submanifold, defined by transversality, of a manifold without boundary. To prove compactness, we must show that $Z_1$ has no limit points where the $x_i$ come together or the $a_{ij}$ tend to zero or infinity. It is easy to eliminate the possibility that the $x_i$ come together by the equivariance of $F$, as $F(x_j,x_k)=-F(x_k,x_j)$ means that the three vectors given by $F$ cannot all point the same way if two of the $x_i$ are the same. The $a_{ij}$ cannot go to infinity since the image of $F:M\times M\setminus\Delta\rightarrow\RR^n$ is bounded by compactness of $M$. The $a_{ij}$ cannot go to zero because $F$ gives a non-zero vector for each pair of distinct points in $M$, and since we have ruled out the possibility that the $x_i$ come together (which is the only way $F$ can be zero), the image of $F$ is bounded away from zero  outside a neighborhood of $\Delta\subset M\times M$. This shows that $Z_1$ is compact. The bundle isomorphism is given to us by transversality:

$$TZ_1\oplus \F^*(\RR^n\times\RR^n)\cong T(M^3\setminus\Delta)\oplus T(\RR^2_{>0}).$$

The isomorphism $T(\RR^2_{>0})\cong \widetilde{P}$ is given by the map $(a,b)\mapsto ae_2+be_3$, and the isomorphism $\F^*(\RR^n\times\RR^n)\cong n\widetilde{P}$ is induced by the map $\RR\times\RR\rightarrow \widetilde{P}$ given by $(a,b)\mapsto ae_1+b(e_2-e_3)$. Both of these isomorphisms are $\Sigma_2^{23}$-equivariant. Hence we have a $\Sigma_2^{23}$-equivariant isomorphism

\begin{equation}\label{zbundleiso}
TZ_1\oplus p^*n\widetilde{P}\rightarrow p^*T(M^3\setminus\Delta)\oplus p^*\widetilde{P}.
\end{equation} 
$\Box$\end{pf}

\begin{defn}\label{Zdef}
Define $Z=(\Sigma_3\times_{\Sigma_2}Z_1)/\Sigma_3$.
\end{defn}





To define $\k$ for a $k$-simplices is straightforward. A $k$-simplex of $\ivmap^{\Sigma_2}(M\times M,\RR^n)$ is an isovariant $F_k:M\times M\times \Delta^{k}\rightarrow\RR^n\times\Delta^{k}$, and the relevant manifold is the $(3m-2n+2+k)$-dimensional submanifold $Z_{1,k}=\F_k^{-1}(0\times \Delta^{k}\subset (M^3\setminus\Delta)\times \RR^2_{>0}\times \Delta^k$. As before, $Z_{1,k}$ is compact, has a $\Sigma_2^{23}$-equivariant proper map to $M^3\setminus\Delta$, and is transverse to $(M^3\setminus\Delta)\times \RR^2_{>0}\times \del_S\Delta^k$ for all $S$. Moreover, there is an isomorphism $TZ_{1,k}\oplus p^*n\widetilde{P}\cong p^*T(M^3\setminus\Delta)\oplus p^*\widetilde{P}\oplus T\Delta^{k}$. We then set $Z_k=(\Sigma_3\times_{\Sigma_2^{23}} Z_{1,k})/\Sigma_3$

This proves

\begin{lem}
There is a well-defined map $\k:\mathcal{T}_2\emb(M^m,\RR^n)\rightarrow C_3(M^m)$ given by $\k(F)=Z$, where $Z=\left(\Sigma_3\times_{\Sigma_2^{23}}\F^{-1}(0\times 0)\right)/\Sigma_3$.
\end{lem}

\section{Proof of Theorem \ref{Muns}}\label{proof}

We now restate the main theorem for convenience and proceed to prove it.

\addtocounter{thm}{-26}

\begin{thm}
The following square is $(3n-4m-5)$-cartesian:

$$\xymatrix{
\emb(M^m,\RR^n) \ar[rr]\ar[dd] & &  \ast \ar[dd] \\
 & & \\
\mathcal{T}_2\emb(M^m,\RR^n)  \ar[rr]^-{\k}& & C_3(M^m)\\
}
$$
\end{thm}

The map $\ast\rightarrow C_3(M)$ assigns to the point $\ast$ the empty manifold. An immediate corollary is

\begin{cor}
Given an isovariant map $F\in \mathcal{T}_2\emb(M^m,\RR^n)$, there exists a manifold

$$Z\in\Om_{3m-2n+2}^{(n-1)P-T{{M}\choose3}}{{M}\choose3}$$

which represents the obstruction to lifting $F$ to $\mathcal{T}_3\emb(M^m,\RR^n)$. If $3n-4m-5\geq 0$, then if $Z$ is null-cobordant there exists an embedding of $M$ in $\RR^n$.
\end{cor}

\addtocounter{thm}{24}

It follows from our Theorem \ref{Muns} that our manifold $Z$ represents the only obstruction to lifting from $F\in\mathcal{T}_2\emb(M,\RR^n)$ to $\mathcal{T}_3\emb(M,\RR^n)$, because there is a $(3n-4m-5)$-connected map $\emb(M,\RR^n)\rightarrow\mathcal{T}_3\emb(M,\RR^n)$. In section \ref{reducetothree} we reduce Theorem \ref{Muns} to the case where $M$ contains exactly three points by an induction argument inspired by the proof of Theorem 5.1 of \cite{we1}, as we have already mentioned. Then in section \ref{specialcase} we prove Lemma \ref{Mthreepoints}, which proves our theorem when $M$ consists of exactly three points.

\subsection{The handle induction}\label{reducetothree}

\begin{pf}[of Theorem \ref{Muns}]
We consider all spaces as images of corresponding cofunctors from $\mathcal{O}(M)$ to Spaces. We will induct on the handle dimension $k$ of open sets $U\in \mathcal{O}(M)$ which are the interior of smooth compact codimension $0$ handlebodies, and finally specialize to $U=M$. Recall that a manifold has handle dimension $k$ if it admits a handle decomposition with handles of at most index $k$. We will prove that if $U$ can be made from handles of index at most $k$, then the square

$$\xymatrix{
\emb(U,\RR^n) \ar[rr]\ar[dd] & &  \ast \ar[dd] \\
 & & \\
\mathcal{T}_2\emb(U,\RR^n)  \ar[rr]& & C(U,\RR^n)\\
}
$$

\begin{figure}
\begin{center}

\input{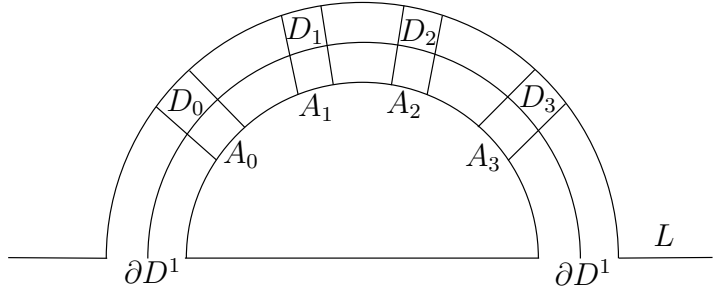}
\caption{The $D_i$ and the $A_i$ for a $1$-handle $D^1\times D^1$ attached along $\del D^1\times D^1$. The $D_i\subset D^1$, and $A_i=D_i\times D^1\subset D^1\times D^1$. Note that removing $k\geq 1$ of the $A_i$ leaves a manifold with $(k-1)$ extra $0$-handles, but one fewer $1$-handle}
\label{handles}
\end{center}
\end{figure}

is $(3n-4k-5)$-cartesian. We will omit the second variable $\RR^n$ from our notation. The base case $k=0$, when $U$ is a tubular neighborhood of a finite set, will established in Lemma \ref{Mthreepoints} below. Let $k>0$ and assume the result for $l<k$. Let $L$ be a smooth compact codimension zero submanifold of $M$, with interior $U$, and let $s>0$ be the number of handles of index $k$. Let $e_j:D^{n-k}\times D^k\rightarrow L$ denote each of the $k$-handles for $j=1$ to $s$. Assume that $e_j^{-1}(\del L)=\del D^{n-k}\times D^k$ for all $j$. Since $k>0$, we may choose, for each $j$, closed pairwise disjoint disks $D_0,D_1,D_2,D_3$ in the interior of $D^k$, and set

$$A^j_i=e_j(D^{n-k}\times D_i)\cap U$$

for each $i$ (see figure \ref{handles}). Then $A^j_i$ is closed in $U$ and if we set $A_i=\cup_j A^j_i$, then $U-A_i$ is the interior of a smooth compact codimension $0$ manifold which admits a handle decomposition with no handles of index greater than or equal to $k$. The same is true for $U_S=\cap_{i\in S}U-A_i$ for each nonempty subset $S\subset\{0,1,2,3\}$. By induction, for each $S\neq \emptyset$ the square

$$\xymatrix{
\emb(U_S) \ar[rr]\ar[dd] & &  \ast \ar[dd] \\
 & & \\
\mathcal{T}_2\emb(U_S)  \ar[rr]& & C(U_S)\\
}
$$

is $(3n-4(k-1)-5)$-cartesian. Hence for each nonempty $S$, the map

$$\emb(U_S)\rightarrow \hofiber(\mathcal{T}_2\emb(U_S)\rightarrow C(U_S))$$

is $(3n-4(k-1)-5)$-connected. Now consider the square diagram

$$\xymatrix{
\emb(U_S) \ar[rr]\ar[dd] & &  \hofiber(\mathcal{T}_2\emb(U_S)\rightarrow C(U_S)) \ar[dd] \\
 & & \\
\holim_{S\neq\emptyset}\emb(U_S)  \ar[rr]& &\holim_{S\neq\emptyset} \hofiber(\mathcal{T}_2\emb(U_S)\rightarrow C(U_S))\\
}
$$

We want to show that the upper horizontal map is $(3n-4k-5)$-connected for all $S$, including $S=\emptyset$. Since $\mathcal{T}_2\emb(-,\RR^n)$ is polynomial of degree $\leq 2$, $C_3(-)$ is polynomial of degree $\leq 3$, and the sets $A_i$ are four in number, the rightmost vertical map is $\infty$-connected. By Goodwillie-Klein \cite{gk}, the leftmost vertical map is $(3n-4k-5)$-connected. By induction and Proposition 1.22 of \cite{tg2}, the lower horizontal map is $(3n-4k-4)$-connected. It follows that the upper horizontal map is $(3n-4k-5)$-connected. Specializing to $U=M$ gives the desired result.


\subsection{Proof of the theorem when $M$ is three points}\label{specialcase}

Now we prove the theorem in the case $k=0$, which is when $U$ is an open tubular neighborhood of a finite set of points. Since $C_3$ is a polynomial of degree $\leq 3$, we can, using the same handle induction argument as above, reduce to the case when $U$ is a tubular neighborhood of at most three points. By Proposition \ref{cobtn}, we may replace the tubular neighborhood $U$ with its zero section $S$. If $S$ has less than three points, then $C(S)$ is contractible, and $\emb(S)\rightarrow\mathcal{T}_2\emb(S)$ is an equivalence. We are thus reduced to proving this theorem in the case where $S=\{x_1,x_2,x_3\}$ consists of exactly three points. For the remainder of this section, $\emb(S)$ will denote $\emb(\{x_1,x_2,x_3\})$ and $C_3(S)$ will denote $C_3(\{x_1,x_2,x_3\})$.
$\Box$\end{pf}

\begin{lem}\label{Mthreepoints}
The square

$$\xymatrix{
\emb(S) \ar[rr]\ar[dd] & &  \ast \ar[dd] \\
 & & \\
\mathcal{T}_2\emb(S)  \ar[rr]& & C_3(S)\\
}
$$

is $(3n-5)$-cartesian, where $S=\{x_1,x_2,x_3\}$.
\end{lem}

This lemma says that the homotopy groups of the homotopy fiber of the left vertical map are isomorphic with the homotopy groups of the right vertical map through a range. The space $\emb(S)$ is the configuration space of three points in $\RR^n$, which has been extensively studied. We will identify $\mathcal{T}_2\emb(S)$ in the next section. The proof of this lemma is broken up into two main steps. In section \ref{hofiber}, we explicitly identify $\hofiber\left(\emb(S)\rightarrow \mathcal{T}_2\emb(S)\right)$, and establish that there is a $(3n-5)$-connected map 

$$S^{2n-3}\rightarrow \hofiber\left(\emb(S)\rightarrow \mathcal{T}_2\emb(S)\right).$$ 

We then identify $\Om C_3(S)$ with $\Om QS^{2n-2}$, and it follows that if the composed map $S^{2n-3}\rightarrow \Om C_3(S)$ induces an isomorphism on $\pi_{2n-3}$, then it is in fact $(4n-5)$-connected. Finally, in section \ref{generator}, we establish the isomorphism between $\pi_{2n-3}\hofiber\left(\emb(S)\rightarrow \mathcal{T}_2\emb(S)\right)$ and $\pi_{2n-3}\Om C_3(S)$ on $\pi_{2n-3}$.

\subsubsection{The homotopy fiber of $\emb(S)\rightarrow \mathcal{T}_2\emb(S)$ and the identification of $\Om C_3(S)$}\label{hofiber}

\begin{lem}\label{hofiberlemma}
For $S=\{x_1,x_2,x_3\}$, there is an equivalence $$\hofiber(\emb(S)\rightarrow\mathcal{T}_2\emb(S))\simeq \hofiber(S^{n-1}\vee S^{n-1}\rightarrow S^{n-1}\times S^{n-1}).$$
\end{lem}

\begin{pf}
Since $\emb(\{x_1,x_2\})\simeq S^{n-1}$, there is a fibration


$$\xymatrix{
S^{n-1}\vee S^{n-1} \ar[rr] & & \emb(\{x_1,x_2,x_3\}) \ar[dd] \\
 & & \\
 & & S^{n-1}\\
}
$$

Note that $\mathcal{T}_2\emb(\{x_1,x_2,x_3\})\simeq S^{n-1}\times S^{n-1}\times S^{n-1}$. Recalling our model for $\mathcal{T}_2\emb(M^m,\RR^n)$, we see that $F$ only needs to specify a non-zero vector of $\RR^n$ for each two element subset of $M$ in an equivariant way. Hence we also have a trivial fibration

$$\xymatrix{
S^{n-1}\times S^{n-1} \ar[rr] & & \mathcal{T}_2\emb(\{x_1,x_2,x_3\}) \ar[dd] \\
 & & \\
 & & S^{n-1}\\
}
$$

The map $\emb(S)\rightarrow\mathcal{T}_2\emb(S)$ induces a map of fibrations

$$\xymatrix{
S^{n-1}\vee S^{n-1} \ar[rr]\ar[dd] & & S^{n-1}\times S^{n-1} \ar[dd] \\
 & & \\
\emb(S)\ar[rr]\ar[dd] & & \mathcal{T}_2\emb(S)\ar[dd]\\
 & & \\
S^{n-1}\ar[rr]_-{=} & & S^{n-1}\\
}
$$

and hence an equivalence

$$\hofiber(\emb(S)\rightarrow\mathcal{T}_2\emb(S))\simeq \hofiber(S^{n-1}\vee S^{n-1}\rightarrow S^{n-1}\times S^{n-1}).$$
$\Box$\end{pf}

There is a homotopy equivalence $\hofiber(S^{n-1}\vee S^{n-1}\rightarrow S^{n-1}\times S^{n-1})\simeq \Sigma(\Om S^{n-1}\wedge\Om S^{n-1})$. One can see this by first identifying the homotopy fiber with $\Om S^{n-1}\ast\Om S^{n-1}$, where $\ast$ denotes the join construction, which in turn maps to $\Sigma(\Om S^{n-1}\wedge\Om S^{n-1})$ by a homotopy equivalence. For details, see \cite{selick}.


The $(2n-3)$-skeleton for the James model for $\Om S^{n-1}$ is $S^{n-2}$, and hence the $(3n-5)$-skeleton of $\Sigma(\Om S^{n-1}\wedge\Om S^{n-1})$ is $S^{2n-3}$. Hence we have a $(3n-5)$-connected map $S^{2n-3}\rightarrow\hofiber(S^{n-1}\vee S^{n-1}\rightarrow S^{n-1}\times S^{n-1})$ given by the inclusion of this skeleton.

%
%
%
%
%
%
%

\begin{lem}
There is an equivalence $\Om C_3(S)\simeq \Om QS^{2n-2}$.
\end{lem}

\begin{pf}
Observe that for $S=\{x_1,x_2,x_3\}$, ${S\choose 3}=\ast$ and by Remark \ref{cobspacethomspace}, $\Om C_3(S)\simeq \Om QS^{2n-2}$ because the tangent bundle to the configuration space is the zero bundle and the bundle $(n-1)P$ is trivial of rank $2n-2$.$\Box$
\end{pf}

By a computation involving the obstruction $Z$, we show that the composed map $S^{2n-3}\rightarrow \Om C_3(S)$ (still to be defined) induces an isomorphism on $\pi_{2n-3}$. But the homology groups of both spaces vanish up to dimension $4n-5$, so using the Hurewicz theomem, the map is actually $(4n-5)$-connected.

\subsubsection{A generator for $\pi_{2n-3}$}\label{generator}

Here we will give a generator of $\pi_{2n-3}\hofiber(\emb(S)\rightarrow\mathcal{T}_2\emb(S))$ and show that the composed map to $\pi_{2n-3}\Om C_3(S)$ generates this group as well. A single point with, of necessity, trivial bundle data will represent a generator of this latter group.

Write $\RR^n=\RR^{n-1}\times\RR$, and let $p_1=(0,1/2)$ and $p_2=(0,-1/2)$ be points in $\RR^n$ in these coordinates. Lemma \ref{hofiberlemma} gives an equivalence of pairs $(S^{n-1}\vee S^{n-1},S^{n-1}\times S^{n-1})\rightarrow(\emb(S),\mathcal{T}_2\emb(S))$. It factors through the inclusion $(S^{n-1}\vee S^{n-1},S^{n-1}\times S^{n-1})\rightarrow (\RR^n\setminus\{p_1,p_2\},\RR^n\setminus\{p_1\}\times\RR^n\setminus\{p_2\})$, where the wedge point is the origin and the spheres are centered around $p_1$ and $p_2$.

Consider the following commutative diagram.


$$\xymatrix{
S^{2n-3}  \ar[rr]^-{\a}\ar[dd]_-{\io} & & \RR^n\setminus\{p_1,p_2\}\ar[dd]_-{(\io_1,\io_2)}\ar[rr] & & \emb(S)\ar[rr]\ar[dd] & & \ast\ar[dd]\\
 & & & & & & \\
D^{2n-2}\ar[rr]_-{\a_{1/2}\times\a_{-1/2}} & & \RR^n\setminus\{p_1\}\times\RR^n\setminus\{p_2\}\ar[rr] & & T_2\emb(S)\ar[rr]_-{\k} & & C_3(S) \\
}
$$


Here $\io$ is the inclusion $S^{2n-3}\rightarrow D^{2n-2}$, and $\io_j$ is the inclusion of $\RR^n\setminus\{p_1,p_2\}$ in $\RR^n\setminus\{p_j\}$ for $j=1,2$.

\begin{defn}\label{alpha}
The map $\a:S^{2n-3}\rightarrow \RR^n\setminus\{p_1,p_2\}$ is given by
$$
\a(v,w) = (\abs{v}^2w + \abs{w}^2v, \abs{v}^2 - \abs{w}^2),
$$
where $(x,y)\in \RR^{n-1} \times \RR$ denotes a point in $\RR^n$, and $S^{2n-3}$ is the unit sphere $\abs{v}^2 + \abs{w}^2 =1$ in $\RR^{n-1} \times \RR^{n-1}$. 
\end{defn}

There are three obvious ways to extend this map over all of $D^{2n-2}$. One is just to extend it by the same formula, which we will also call $\a$. The other two make use of the fact that $\abs{v}^2 + \abs{w}^2 = 1$ on the sphere, so that we may write the restriction of $\a$ to the sphere in two equivalent ways there, and extend them over the whole disk in the obvious way. These maps are denoted $\a_{1/2}$ and $\a_{-1/2}$, and are given by

$$
\a_{1/2}(v,w) = (w+\abs{w}^2(v-w), 1 - 2\abs{w}^2)
$$

and

$$
\a_{-1/2}(v,w) = (v+\abs{v}^2(w-v), 2\abs{v}^2 - 1).
$$

The following lemma verifies these maps have the target we claim they do.

\begin{lem}
The restriction of $\a$ to the sphere misses $p_1=(0,1/2)$ and $p_2=(0,-1/2)$, and the map $\a_{1/2}$ on the whole disk misses the point $(0,1/2)$, and likewise for $\a_{-1/2}$ and the point $(0,-1/2)$.
\end{lem}

\begin{pf}
For the first fact, note that there are only two ways in which $\abs{v}^2w+\abs{w}^2v$ can be zero: one of the coordinates $v$ or $w$ is zero, or $\abs{v}=\abs{w}$. In the first case, the other coordinate must have length one, in which case $\abs{v}^2-\abs{w}^2$ is $\pm 1$, and in the second, $\abs{v}^2-\abs{w}^2=0$. For the second fact, the argument is the same for both $\a_{1/2}$ and $\a_{-1/2}$, so we will argue only that $\a_{1/2}$ misses $(0, 1/2)$. If $\a_{1/2}(v,w)=(0,1/2)$, then we must have $\abs{w}^2=1/4$. Solving for $v$ in terms of $w$ we obtain $v=(1-1/\abs{w}^2)w$ using the first part of the map, and using our previous observation and taking lengths we obtain $\abs{v}=3/2$, which is impossible on $\abs{v}^2+\abs{w}^2\leq1$.
$\Box$
\end{pf}

\begin{lem}
The map of pairs $$(\a,\a_{1/2}\times\a_{-1/2}):(S^{2n-3},D^{2n-2})\rightarrow(\RR^n\setminus\{p_1,p_2\},\RR^n\setminus\{p_1\}\times\RR^n\setminus\{p_2\})$$ represents a generator of $\pi_{2n-3}\hofiber(S^{n-1}\vee S^{n-1}\rightarrow S^{n-1}\times S^{n-1})$.
\end{lem}

\begin{pf}
It is known (see, for example, \cite{bcss}) that the map assigning to each smooth map $f:S^{2n-3}\rightarrow S^{n-1}\vee S^{n-1}$ the linking number $lk(f^{-1}(y_1),f^{-1}(y_2))$, where $y_1\in S^{n-1}\vee \ast$ and $y_2\in \ast\vee S^{n-1}$ are regular values of $f$, provides an isomorphism of $\pi_{2n-3}\hofiber(S^{n-1}\vee S^{n-1}\rightarrow S^{n-1}\times S^{n-1})$ with $\ZZ$. The points $(0,1)$ and $(0,-1)$ are regular values of $\a$. The inverse images of both points are $(n-2)$-dimensional spheres $S_{\pm1}^{n-2}=\a^{-1}(0,\pm1)$. One easily sees that $\a^{-1}(0,1)=\{\abs{v}^2=1\}$ and $\a^{-1}(0,-1)=\{\abs{w}^2=1\}$. The linking number of these spheres is $1$. This can be computed by counting intersections of one of the spheres with a bounding disk. If we let $D_{+1}^{n-1}=\{\abs{v}^2\leq1\}$, then $\del D_{+1}^{n-1}=S_{+1}^{n-2}$, and this disk intersects $S_{-1}^{n-2}$ only at $(v,w)=(0,0)$.
$\Box$\end{pf}

We now explicitly construct the manifold $Z$. Recall that $Z$ is constructed by determining when the three vectors determined by evaluating $F$ on pairs of a triple $(x_1,x_2,x_3)$ of distinct points of $M$ point in the same direction. We have a parametrized family of maps $F_s$, parametrized by coordinates $s=(v,w)$ in the disk $D^{2n-2}$. The maps $F_s$ are easy to describe, since $M=\{x_1,x_2,x_3\}$ contains just three points and we have explicitly described the map $D^{2n-2}\rightarrow \RR^{n}\setminus\{p_1\}\times\RR^{n}\setminus\{p_2\}$.

\begin{lem}
The equations 

\begin{eqnarray*}
F_s(x_1,x_2) &=& \a_{1/2}(v,w)-(0,1/2)\\
F_s(x_2,x_3) &=& (0,1)\\
F_s(x_3,x_1) &=& (0,-1/2)-\a_{-1/2}(v,w)
\end{eqnarray*}

represent the composed map $D^{2n-2}\rightarrow \mathcal{T}_2\emb(S)$ and define nonzero vectors for each $s=(v,w)\in D^{2n-2}$. Moreover, the map $\F_s:\RR^2_{>0}\times D^{2n-2}\rightarrow\RR^n\times\RR^n$, a parametrized family of maps defined by the above using Definition \ref{Fdef}, is transverse to $0\times 0$, and its only zero occurs when $(v,w)=(0,0)$.
\end{lem}

\begin{pf}
The properties of $\a_{1/2}$ and $\a_{-1/2}$ noted above ensure that this defines a non-zero vector for each $(s,x_i, x_j)$ for $i\neq j$. To find the zeroes of $\F_s$, that is, to compute the manifold $Z$, we need to compute when the $F_s(x_i,x_j)$ are positive multiples of $(0,1)$. This is the case if $(v,w)=(0,0)$. We claim that this is the only solution.

If $v=0$, then $\abs{w}=1$ since the first coordinate of $F_{(0,w)}(x_1,x_2)$ must be zero, but in this case the second coordinate is negative. By symmetry this rules out the possibility that there is a solution when either $v=0$ or $w=0$. Now assume that $v,w\neq 0$. Again considering that the first coordinate of $F_{(v,w)}(x_1,x_2)$ must be zero, we see that there must be a linear dependence between $v$ and $w$. In particular, we must have $v=(1-1/\abs{w}^2)w$, and $w=(1-1/\abs{v}^2)v$. By substitution and algebra we end up seeking solutions to $2\abs{w}^4-3\abs{w}^2+1=0$, which are $\abs{w}^2=1$ or $\abs{w}^2=1/2$. When $\abs{w}^2=1$, we must have $v=0$, which has already been ruled out. When $\abs{w}^2=1/2$, the second coordinate of $F_s(x_1,x_2)$ is negative. Hence $v=w=0$ is the only solution.

To check that $\F_s$ is transverse to $0\times 0$ amounts to checking that the matrix $D\F_s$ has rank $2n$.

Write $\F_s=((y_1,u_1),(y_2,u_2))$. Then

$$D\F_s=
\left(
\begin{array}{cccc}
\frac{dy_1}{dv} & \frac{dy_1}{dw} & \frac{dy_1}{da_{31}} & \frac{dy_1}{da_{12}}\\
\frac{du_1}{dv} & \frac{du_1}{dw} & \frac{du_1}{da_{31}} & \frac{du_1}{da_{12}}\\
\frac{dy_2}{dv} & \frac{dy_2}{dw} & \frac{dy_2}{da_{31}} & \frac{dy_2}{da_{12}}\\
\frac{du_2}{dv} & \frac{du_2}{dw} & \frac{du_2}{da_{31}} & \frac{du_2}{da_{12}}\\
\end{array}
\right)
$$

Since $\F_s$ and $\F_s'$ are related by an invertible linear transformation, it is enough to check that $\F_s'$ has rank $2n$. Letting $I_k$ denote the $k\times k$ identity matrix, we find that

$$D\F'|_{(0,0)}=
\left(
\begin{array}{cccc}
a_{31}I_{n-1} & 0 & 0 & 0 \\
0 & 0 & 1/2 & 0 \\
0 & a_{12}I_{n-1} & 0 & 0 \\
0 & 0 & 0 & 1/2
\end{array}
\right)
$$

which has rank $2n$, since the $a_{ij}$ of Definition \ref{aij} are positive.
$\Box$
\end{pf}

This completes the proof of Lemma \ref{Mthreepoints}, as we have shown that a generator of $\pi_{2n-3}\left(\hofiber(S^{n-1}\vee S^{n-1}\rightarrow S^{n-1}\times S^{n-1})\right)$ goes to a generator of the cobordism group $\Om_{3m-2n+2}^{(n-1)P-T{M\choose3}}{M\choose3}$ by this construction.

\section{Smooth knotting of spheres}\label{smoothknot}

As an application of our Theorem \ref{Muns}, we recover results due to Haefliger in \cite{haef2} on the knotting of smooth spheres. We should note, however, that he used surgery theory to prove these, and it was important that the manifolds to be embedded were spheres. As an application of our Theorem \ref{Muns}, we will prove

\begin{thm}[\cite{haef2}, 8.14]\label{smoothknotthm}
$\pi_0\emb(S^{2k+1},\RR^{3k+3})$ is isomorphic with $\ZZ$ if $k$ is odd, and $\ZZ/2$ if $k$ is even.
\end{thm}

Dax's Theorem \ref{Dax} says that the map $\emb(M^m,\RR^n)\rightarrow\mathcal{T}_2\emb(M^m,\RR^n)$ is $(2n-3m-3)$-connected. Kervaire \cite{ker} proves  that $\pi_0\mathcal{T}_2\emb(S^m,\RR^n)=0$ for $2n-3m-1>0$, and hence all embeddings of $S^m$ in $\RR^n$ are isotopic if $2n-3m-3>0$. This will play an important role in enumerating embeddings of $S^m$ in $\RR^n$ when $2n-3m-3=0$. To prove Theorem \ref{smoothknotthm}, we also need to know about $\pi_1\mathcal{T}_2\emb(S^m,\RR^n)$.

\begin{lem}\label{t2lemma}
$\pi_k\mathcal{T}_2\emb(S^{2k+1},\RR^{3k+3})=0$ for $k=0,1$.
\end{lem}

The proof of Lemma \ref{t2lemma} will occupy most of the rest of this section. Theorem \ref{smoothknotthm} follows easily from this lemma and our Theorem \ref{Muns}.
 
Denote by $\Gamma(X,Y)$ the space of sections of some fibration over $X$ with fibers $Y$. The space $\mathcal{T}_1\emb(S^m,\RR^n)$ is weakly equivalent to $\Gamma(S^m,V_{m,n})$ by the Smale-Hirsch theorem, where $V_{m,n}$ is the Stiefel manifold of $m$-frames in $\RR^n$. The fibration in question has as its total space the space of vector bundle monomorphisms from $TS^m$ to $T\RR^n$. Recall that the map $\mathcal{T}_2\emb(S^m,\RR^n)\rightarrow\Gamma(S^m,V_{m,n})$ restricts an isovariant map $F$ to the diagonal and records the induced map of normal bundles. Let $X$ and $Y$ be spaces with a $\Sigma_2$ action. Denote by $\equ(X,Y)$ the space of $\Sigma_2$-equivariant maps from $X$ to $Y$. Consider the space $\equ(S^m\times S^m\setminus\Delta, S^{n-1})$, where $\Sigma_2$ acts by switching the coordinates in the first variable and antipodally in the second.  We can restrict an equivariant map $S^m\times S^m\setminus \Delta\rightarrow S^{n-1}$ to the bundle of $(m-1)$ spheres associated to a tubular neighborhood of the diagonal $\Delta\subset S^m\times S^m$. Since $\Delta\cong S^m$, we can view this as giving an equivariant map $S^{m-1}\rightarrow S^{n-1}$ for each point in the diagonal. This can be interpreted as a section of a bundle over $S^m$ whose fibers are $\equ(S^{m-1},S^{n-1})$. It is built in exactly the same way the bundle of vector bundle monomorphisms of $TS^m$ in $T\RR^n$ is built from $S^m$ and $V_{m,n}$; we replace $V_{m,n}$ with $\equ(S^{m-1},S^{n-1})$.

There is a map $\mathcal{T}_2\emb(S^m,\RR^n)\rightarrow\equ(S^m\times S^m\setminus\Delta, S^{n-1})$ given by sending an isovariant map $F=(f_1,f_2):S^m\times S^m\rightarrow \RR^n\times\RR^n$ to the restriction of $\frac{f_1-f_2}{\abs{f_1-f_2}}$ to the complement of the diagonal $\Delta$. Likewise, there is a map $\Gamma(S^m,V_{m,n})\rightarrow\Gamma(S^m,\equ(S^{m-1},S^{n-1}))$ induced by the inclusion $V_{m,n}\rightarrow\equ(S^{m-1},S^{n-1})$ which associates a linear length preserving map of rank $m$ to an equivariant map of spheres (with antipodal actions) by restriction.



\begin{lem}\label{equmodel}
The square diagram
$$\xymatrix{
\mathcal{T}_2\emb(S^m,\RR^n)  \ar[rr]\ar[dd] & & \equ(S^m\times S^m\setminus\Delta, S^{n-1}) \ar[dd] \\
 & & \\
\Gamma(S^m,V_{m,n})  \ar[rr] & & \Gamma(S^m,\equ(S^{m-1},S^{n-1}))\\
}
$$
 is homotopy cartesian.
\end{lem}

\begin{pf}
By Theorem 9.2 of \cite{we1}, the left vertical fibers are equivalent to $\equ_c(S^m\times S^m\setminus\Delta,S^{n-1})$, where the subscript $c$ denotes the additional requirement that the sections should be given in a neighborhood of the diagonal $\Delta$. By inspection, this is the right vertical fiber.
\end{pf}

If follows from Lemma \ref{equmodel} that the connectivity of the top vertical map is the same as that of the bottom vertical map.

\begin{lem}[\cite{haefh2}, Lemma 1.1]\label{stiefelrange}
The map $V_{m,n}\rightarrow \equ(S^{m-1},S^{n-1})$ is $(2n-2m-1)$-connected.
\end{lem}




%



It follows from Lemma \ref{stiefelrange} that

\begin{thm}[\cite{haefh2}, Theorem 4.2]
The map $$\Gamma(S^m,V_{m,n})\rightarrow\Gamma(S^m,\equ(S^{m-1},S^{n-1}))$$ is $(2n-3m-1)$-connected.
\end{thm}

This follows from the fact that if $E\rightarrow B$ is a fibration with $k$-connected fiber and B a $d$-dimensional CW-complex, then the space of sections is $(k-d)$-connected. Hence we have proven

\begin{thm}
The map $$\mathcal{T}_2\emb(S^m,\RR^n)\rightarrow\equ(S^m\times S^m - \Delta, S^{n-1})$$ is $(2n-3m-1)$-connected.
\end{thm}






We may replace $\equ(S^m\times S^m\setminus\Delta,S^{n-1})$ with $\equ(S^m,S^{n-1})$ because the map from $S^m\times S^m\setminus\Delta\rightarrow S^m$ which sends $(x,y)\rightarrow \frac{x-y}{\abs{x-y}}$ is an equivariant homotopy equivalence, with homotopy inverse $x\mapsto (x,-x)$. By Lemma \ref{stiefelrange}, we have a $(2n-2m-3)$-connected map $V_{m+1,n}\rightarrow\equ(S^m,S^{n-1})$, and $V_{m+1,n}$ itself is $(n-m-2)$-connected.

Now let $m=2k+1, n=3k+3$. The map $$\eta_2:\emb(S^{2k+1},\RR^{3k+3})\rightarrow\mathcal{T}_2\emb(S^{2k+1},\RR^{3k+3})$$ is $0$-connected, meaning it is surjective on components, but the map 
$$\emb(S^{2k+1},\RR^{3k+3})\rightarrow\mathcal{T}_3\emb(S^{2k+1},\RR^{3k+3})$$ is $k$-connected, and hence gives an isomorphism on $\pi_0$ when $k\geq 1$.
The map

$$\mathcal{T}_2\emb(S^{2k+1},\RR^{3k+3})\rightarrow\equ(S^{2k+1}, S^{3k+2})$$

is $2$-connected, $V_{2k+2,3k+3}\rightarrow\equ(S^{2k+1}, S^{3k+2})$ is $(2k+3)$-connected, and $V_{2k+2,3k+3}$ is itself $k$-connected. It follows that 

$$\pi_0\mathcal{T}_2\emb(S^{2k+1},\RR^{3k+3})=\pi_1\mathcal{T}_2\emb(S^{2k+1},\RR^{3k+3})=0.$$ This completes the proof of Lemma \ref{t2lemma}.

Now consider the long exact sequence of homotopy groups of the fibration $L_3' \rightarrow \emb \rightarrow \mathcal{T}_2\emb$. From our Theorem \ref{Muns} we have a $k$-cartesian square

$$\xymatrix{
\emb(S^{2k+1},\RR^{3k+3}) \ar[rr]\ar[dd] & & \ast \ar[dd] \\
 & & \\
\mathcal{T}_2\emb(S^{2k+1},\RR^{3k+3}) \ar[rr] & &  C(S^{2k+1},\RR^{3k+3})\\
}
$$


Hence there is a $k$-connected map of vertical fibers $L_3' \rightarrow \Om C$. Taking $\pi_0$, we see that $\pi_0 \Om C= \Om^{(3k+2)P-T{S^{2k+1}\choose 3}}_0{S^{2k+1}\choose 3}$. If $k$ is odd, then this group is $\ZZ$, and when $k$ is even it is $\ZZ/2$.

To explain this computation, we need to consider the action of $\pi_1{S^{2k+1}\choose 3}$. This group is isomorphic with $\Sigma_3$ since $S^{2k+1}$ is simply connected and $k\geq1$. Recall that we made $P$ from a representation of $\Sigma_3$, so the homomorphism $w(P):\Sigma_3\rightarrow \{+1,-1\}$ factors through $GL_2(\RR)$ as $\Sigma_3\rightarrow GL_2(\RR)\rightarrow \{+1,-1\}$, where the first map is the representation in question, and the second map records the sign of the determinant. Since the elements of order two generate the group, it is enough to understand $w(P)$ on such elements. Each element $\sigma$ of order two acts by a reflection on the plane, and hence $w(P)(\sigma)=-1$. More generally, $w((3k+2)P)\rightarrow GL_{6k+4}(\RR)$, and for an element $\sigma$ of order two, if $k$ is even, then $w((3k+2)P)(\sigma)=+1$, and if $k$ is odd, then $w((3k+2)P)(\sigma)=-1$. As for the map $w(T{S^{2k+1}\choose 3})\rightarrow\{+1,-1\}$, note that $S^{2k+1}$ is orientable, and hence so is $(S^{2k+1})^3\setminus\Delta$. But any element $\sigma\in\Sigma_3$ of order two changes the sign of the orientation class of $(S^{2k+1})^3\setminus\Delta$ because $2k+1$ is odd. Hence $w(T{S^{2k+1}\choose 3})(\sigma)=-1$ for any element $\sigma$ of order two. It follows that $w((3k+2)P)=w(T{S^{2k+1}\choose 3})$ if $k$ is odd, and $w((3k+2)P)\neq w(T{S^{2k+1}\choose 3})$ if $k$ is even. Proposition \ref{countclasses} implies that $\pi_0 \Om C= \Om^{(3k+2)P-T{S^{2k+1}\choose 3}}_0{S^{2k+1}\choose 3}$ is isomorphic with $\ZZ$ if $k$ is odd, and $\ZZ/2$ if $k$ is even. This completes the proof of Theorem \ref{smoothknotthm}.

\subsection{Acknowledgments}

This paper represents my dissertation completed at Brown University under the guidance of Tom Goodwillie. I would like to thank my readers, Kiyoshi Igusa and George Daskalopoulos, for their helpful comments. Thanks also to Nick Kuhn and Rainer Vogt for pointing me to useful references. I would like to thank the referee for extensive, detailed and helpful comments which greatly improved the exposition of this paper. I am particularly indebted to my advisor, Tom Goodwillie, who has been and continues to be extremely generous and patient in discussing his ideas.

\end{document}